\documentclass[11pt]{amsart}
\usepackage{amscd,amssymb}
\usepackage[dvips]{graphics} 
\usepackage{tabularx}
\input diagrams
\diagramstyle[PostScript=dvips]

\newcommand{\includefig}[1]{\raisebox{-3ex}{\resizebox{!}{7ex}{\includegraphics{pix/#1.eps}}}}

\newcommand{\ds}{\displaystyle}

\newcommand{\tcongi}{\includefig{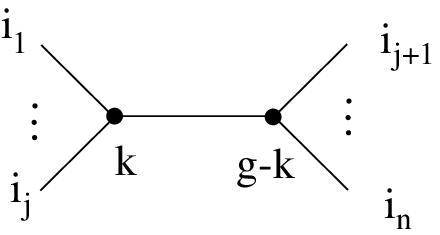}}
\newcommand{\ocongi}{\includefig{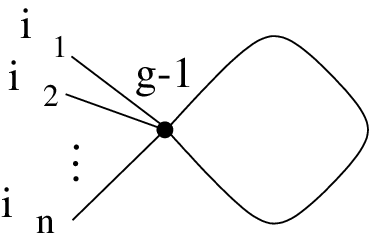}}

\newcommand{\scheme}[1]{\mathcal {#1}}

\newcommand{\mbar}{\M}

\newcommand{\mgnbar}{\overline{\scheme{M}}_{g,n}}              
\newcommand{\mgnrmbar} {\mgnbar^{1/r,\bm}}    
\newcommand{\mgnrbar} {\mgnbar^{1/r}}    
  
\newcommand{\jac}{\text{Jac\,}}        

\newcommand{\tensor}{\otimes}
\newcommand{\cross}{\times}

\newcommand{\irightarrow}{\rTo^{\sim}}

\newcommand{\ck}{{\mathcal K}}
\newcommand{\cl}{{\mathcal L}}
\newcommand{\ce}{{\mathcal E}}
\newcommand{\cf}{{\mathcal F}}
\newcommand{\co}{{\mathcal O}}

\newcommand{\ch}{{\mathcal H}}
\newcommand{\cd}{{\mathcal D}}

\newcommand{\sheafhom}{\ch\kern-.15em om}

\newcommand{\Aut}{\operatorname{Aut}} 
\newcommand{\bgamma}{\boldsymbol{\gamma}}
\newcommand{\bm}{\mathbf{m}}          %
\newcommand{\bbm}{\bar{\bm}}

\newcommand{\bo}{\boldsymbol{1}}

\newcommand{\br}{\mathbf{r}}
\newcommand{\bt}{\mathbf{t}}
\newcommand{\btau}{\boldsymbol{\tau}}

\newcommand{\bx}{\mathbf{x}}          
\newcommand{\cft}{CohFT}              
\newcommand{\cfts}{CohFTs}            
\newcommand{\chr}{\ch^{(r)}}          
\newcommand{\cht}{\tilde{\ch}^{(r)}}          
\newcommand{\cv}{c^{1/r}}            
\newcommand{\ev}{\mathrm{ev}}         

\newcommand{\KdV}{\mathrm{KdV}}
\newcommand{\M}{\overline{\MM}}       
\newcommand{\MM}{\scheme{M}}          

\newcommand{\nc}{{\mathbb{C}}}        
\newcommand{\nq}{{\mathbb{Q}}}        
\newcommand{\nz}{{\mathbb{Z}}}        
\newcommand{\Phit}{\widetilde{\Phi}}  

\def\lift#1#2{
  \dimen0 = \unitlength 
  \multiply\dimen0 by #1 \divide \dimen0 by 2
  \dimen1 = \dimen0 
  \multiply \dimen1 by 7 \divide \dimen1 by 10
  \raise\dimen1
     \hbox{\hskip 0.3cm ${\vbox to \dimen0{}}$ \enspace #2}}

\newtheorem{thm}{Theorem}[section]
\newtheorem{lm}[thm]{Lemma}

\theoremstyle{definition}

\newtheorem{rem}[thm]{Remark}
\newtheorem{rems}[thm]{Remarks}

\newtheorem{df}[thm]{Definition}
\newtheorem{ex}[thm]{Example}

\theoremstyle{remark}

\begin{document}
\addtocounter{section}{-1}

\title[Tensor Products of Frobenius Manifolds]
{Tensor Products of Frobenius Manifolds and Moduli Spaces of Higher Spin
Curves}

\author
[T. J. Jarvis]{Tyler J. Jarvis}
\address
{Department of Mathematics, Brigham Young University, Provo, UT 84602, USA}
\email{jarvis@math.byu.edu}
\thanks{Research of the first author was partially supported by NSA 
grant MDA904-99-1-0039}

\author
[T. Kimura]{Takashi Kimura}
\address
{Department of Mathematics, 111 Cummington Street, Boston
University, Boston, MA 02215, USA} 
\email{kimura@math.bu.edu}
\thanks{Research of the second author was partially supported by NSF grant
  number DMS-9803427}

\author
[A. Vaintrob]{Arkady Vaintrob\vskip 0.2cm
\centerline{\emph{In memory of Mosh\'e Flato}}
}
\address
{Department of Mathematical Sciences, New 
Mexico State University, Las Cruces, NM 88003, USA}
\email{vaintrob@math.nmsu.edu}

\date{\today}

\begin{abstract} 
  
  We review progress on the generalized Witten conjecture and some of
  its major ingredients.  This conjecture states that certain
  intersection numbers on the moduli space of higher spin curves
  assemble into the logarithm of the $\tau$ function of a
  semiclassical limit of the $r$-th Gelfand-Dickey (or $\KdV_r$)
  hierarchy. Additionally, we prove that tensor products of the
  Frobenius manifolds associated to such hierarchies admit a geometric
  interpretation in terms of moduli spaces of higher spin structures.
  We also elaborate upon the analogy to Gromov-Witten invariants of
  a smooth, projective variety.
\end{abstract}
\maketitle

\section{Introduction}
\label{intro}

In recent years, there has been a great deal of interaction between
mathematics and quantum field theory. One such area has been in the
context of topological gravity coupled to topological matter.

The notion of a cohomological field theory (\cft), due to Kontsevich
and Manin \cite{KM1} is an axiomatization of the expected
factorization properties of the correlators appearing in such a
theory, {c.f.} \cite{W,W2,W3,Di}. The Gromov-Witten invariants of a
smooth, projective variety $V$ provide a rigorous construction of a
\cft, one for each $V$, with a state space $H^\bullet(V)$ and a
(Poincar\'e) metric $\eta$. The genus zero correlators endow
$H^\bullet(V)$ with the structure of a (formal) Frobenius manifold
which is a deformation of the cup product, and which yields the
structure of quantum cohomology. Furthermore, for two smooth,
projective varieties $V$ and $V'$, the quantum cohomology of
$H^\bullet(V\,\times\,V')$ is related to
$H^\bullet(V)\,\otimes\,H^\bullet(V')$ through a deformation of the
K\"unneth formula.  The Gromov-Witten invariants are known to be
symplectic invariants of $V$ and are therefore of great mathematical
interest.

When $V$ is a point, Witten conjectured \cite{W3} and Kontsevich proved
\cite{Ko} that intersection numbers of tautological cohomological classes on
the moduli space of stable curves assemble into a $\tau$ function of the
$\KdV$ hierarchy. However, the $\KdV$ hierarchy is only the first of a series
of integrable hierarchies, one for each integer $r\,\geq\,2$, called the
$\KdV_r$ or the $r$-th Gelfand-Dickey hierarchy, where the usual $\KdV$
hierarchy is $\KdV_2$. Witten stated a generalization \cite{W} of this
conjecture where he suggested a construction of moduli spaces and cohomology
classes on them whose intersection numbers assemble into a $\tau$ function of
the $\KdV_r$ hierarchy. At the time that his conjecture was formulated,
however, the relevant compact moduli spaces had not yet been constructed.

In \cite{J,J2}, the moduli space of stable $r$-spin curves of 
genus $g$, $\M_{g,n}^{1/r}$, was constructed. In \cite{JKV}, it 
was proved that certain intersection numbers on $\M_{0,n}^{1/r}$ 
assemble into the genus zero part of the $\tau$ function of the 
$\KdV_r$ hierarchy. In addition, there is a (genus zero part of 
a) \cft\ associated to $\KdV_r$ whose associated Frobenius 
manifold has a potential that is polynomial in its flat 
coordinates. The relevant potential function in genus zero 
satisfies the string equation, dilaton equation, and topological 
recursion relations. These proofs are purely algebro-geometric. 
Axioms have also been formulated for a 
\textsl{virtual class} in higher genus.  The existence of a class
meeting these axioms gives a \cft\ in all genera.  For the case 
of $g=0$ and all $r>1$ as well as for the case of $r=2$ and all 
$g \ge 0$, this class was constructed in \cite{JKV}, but an 
algebro-geometric construction of this class for all $r$ and $g$ 
has yet to be completed. 

From a physical perspective, it should not be surprising that the
theory of Gromov-Witten invariants and the $\KdV_r$ theory should
share many features (see the chart on next page).  Gromov-Witten
invariants correspond to a theory of topological gravity with a 
matter sector arising from the topological sigma model with 
target variety $V$. However, one may also consider other choices 
of matter sector. Witten was led to his $\KdV_r$ conjecture by 
studying the theory of topological gravity coupled to a matter 
sector arising from a coset model. Even mirror symmetry has an 
analog in this theory; namely, the Frobenius manifold given by 
the \cft\ associated to $\KdV_r$ is isomorphic to that arising 
from the  $A_{r-1}$ singularity, \textsl{c.f.}\ 
\cite{DiVV,Du,Ma,Ma2}. 

In this paper, we elaborate upon the analogy between the theory 
of Gromov-Witten invariants and the $\KdV_r$ theory (see the 
chart on the next page) using the results in \cite{JKV}.  We 
construct a geometric ``A-model'' realization (extending the  
results in \cite{JKV}) of the tensor product of the Frobenius 
manifolds associated to $\KdV_r$ in terms of the moduli space of 
multiple spin structures on stable curves. 

\vskip 0.2cm
\textsl{Acknowledgments.} We would like to thank A. Kabanov for useful discussions.

\section{Comparison of \cft s related to Stable Maps and Higher Spin 
  Curves}

This table illustrates some of the similarities and differences
between stable maps and higher spin curves.  
Details and notation will be explained in subsequent sections.

\hspace{-.5in}\begin{tabular}{|m{1.25in}|m{2.375in}|m{2.375in}|}
\hline
&{\bf Stable maps} & {\bf Higher spin curves} \\ \hline \hline
Key Moduli Spaces & Stable maps:  $\M_{g,n}(V)$ & $r$-spin
curves: $\M_{g,n}^{1/r}$ \\ \hline 
State Space & $H^\bullet(V)$ & $\chr$ with basis $\{e_0,\ldots,e_{r-2}\}$
\\ \hline 
Metric & $\eta(\gamma',\gamma'')=\int_V \gamma'\cup \gamma''$ &
$\eta(e_{m'}, e_{m''})=\delta_{m'+m'',r-2}$ \\ \hline 
Flat Identity & $\mathbf{1}$ & $e_0$ \\ \hline
Grading & Integral & Fractional \\ \hline
Fundamental Class & Virtual fundamental class $[ \M_{g,n}(V) ]^\mathrm{virt}$
& Orbifold fundamental class $[ \M_{g,n}^{1/r} ]$ \\ \hline
``Gromov-Witten'' Classes &  True Gromov-Witten classes: $\ev^*(\bgamma) :=
\ev_1^* \gamma_1\cup\ldots\cup \ev_n^* \gamma_n$ & Virtual class $\cv(\bm)$:
constructed for $r=2$ and all $g$ as well as for  $g =0$ and all $r$.\\ \hline  
Gravitational Correlators &
$\left<\tau_{a_1}(\gamma_1)\ldots\tau_{a_n}(\gamma_n)\right>_g$ 
\mbox{$= \int_{[\M_{g,n}(V)]^\mathrm{virt}}
 \ev^*(\bgamma) \prod_{i=1}^n \psi_i^{a_i} $}
&
$\left<\tau_{a_1}(e_{m_1})\ldots\tau_{a_n}(e_{m_n})\right>_g$
\mbox{$=\int_{[\M_{g,n}^{1/r,\bm}]} \cv \prod_{i=1}^n\psi_i^{a_i}$ }
\\ \hline
Free Energy & Large phase space potential:
$\left<\exp(\bt\cdot\btau) \right>$ & Logarithm of tau function:
$\left<\exp(\bt\cdot\btau) \right>$ \\ \hline 
Virasoro Algebra Action & The Virasoro Conjecture: proved for $V$ a point. Proved for
$g=0$ and in some cases for $g=1$.&  $W_r$-algebra conjecture: proved
for $g=0$ as well as for $r=2$ \\ \hline 
Integrable System & $\KdV_2$ when $V$ is a point. Known in $g=0$ and
in some cases in $g=1$. & $\KdV_r$
hierarchy \\ \hline (Generalized) Witten Conjectures & Proved for $V$ a
point. 
& Proved for $g=0$ as well as for $r=2$.
\\ \hline 
Cohomological Field Theory &
$\Lambda_{g,n}=
\mathrm{st}_*(\ev^*(\bgamma)\cap [\M_{g,n}]^\mathrm{virt})$ &
$\Lambda_{g,n}=r^{1-g} p_*\cv(\bm)$ \\ \hline 
WDVV equation & Holds & Holds \\ \hline
String equation & Holds & Holds \\ \hline
TRR & Holds & Holds \\ \hline
Tensor Product of \cft s & Realized by the product $V_1 \cross V_2$
(Quantum K\"unneth theorem in 
$g=0$) & Realized by the fibered
product $\M^{1/r,1/s}_{g,n}$ \\ \hline
Isomorphisms & Mirror symmetry & Versal deformation of $A_{r-1}$
singularity\\ \hline 

\end{tabular}

\begin{rems} \ 

\begin{enumerate}
\item From a geometric point of view the role of the virtual
  fundamental class of stable maps is played more by the virtual class
  $\cv$ than by the orbifold fundamental class.  Indeed, most of the
  real geometric subtlety of the r-spin theory is contained in
  the construction of $\cv$, whereas the smoothness of the moduli
  spaces $\mgnrbar$ implies that the orbifold class is a simple
  modification of the standard fundamental class.

\item For general $V$ the generalized Witten conjecture is not
  precisely formulated, although it is expected that for each $V$ the
  exponential of the free energy is a tau function for some integrable
system.  In $g=0$ and for some cases in $g=1$ this integrable system is
described in \cite{Du2,DuZh}.

\item The $W_r$-algebra contains the Virasoro algebra, and so the
  corresponding conjecture is a generalization of the Virasoro
  Conjecture.

\end{enumerate}
\end{rems}

\section{The $\KdV_r$ hierarchy and its semiclassical limit}

\subsection{The $\KdV_r$ hierarchy}
\label{ssec:kdv}

In this section we introduce the generalized $\KdV$ hierarchies (also called
the Gelfand-Dickey hierarchies) and some of their solutions which will be
used in the subsequent sections.

Fix an integer $r\ge 2$ and consider differential operators of order $r$ in
$D=\frac{i}{\sqrt{r}} \frac{\partial}{\partial x}$ (where the factor
$\frac{i}{\sqrt{r}}$ is added for convenience). A monic differential
operator of order $r$, after conjugation by the operator of  multiplication
by an appropriate function, can be written in the form 
\begin{equation}
  \label{eq:dops}
  L=D^r-\sum_{m=0}^{r-2}u_m(x)D^m,
\end{equation}
where $u_m$ are formal functional variables.
Denote by $\cd$ the affine space of the operators of type~(\ref{eq:dops}).

For every operator $L\in \cd$ there exists a unique pseudo-differential
operator
$$
L^{1/r}=D+\sum_{m>0}w_m D^{-m},
$$
such that $(L^{1/r})^r=L$.
All coefficients $w_m$ of $L^{1/r}$ are differential polynomials in
$u_0,u_1,\ldots,u_{r-2}$. 

A pseudodifferential operator $\ds Q=\sum_{m\ge -n}v_mD^{-m}$
can be decomposed as $Q=Q_+ + Q_-$, where $\ds Q_+=\sum_{m=-n}^0 v_m D^{-m}$
is the differential part of $Q$.

For every $k\ge 0$, the operators $Q=(L^{1/r})^k$ and  $L=(L^{1/r})^r$
commute yielding $[Q_+,L]=-[Q_-,L]$ and, therefore, $[Q_+,L]$ is a
differential operator of order not exceeding $r-2$. Thus the expression 
$ \frac{\partial L}{\partial t} = [Q_+,L]$ 
gives a meaningful system of differential equations for unknown functions
$u_0,\ldots,u_{r-2}$. It can be shown that for different values of $k$ these
equations are consistent (i.e. the corresponding flows on $\cd$ commute) which
allows one to define the $\KdV_r$ (or the $r$-th Gelfand-Dickey) hierarchy as
the following infinite family of differential equations on $\cd$:
\begin{equation}
  \label{eq:kdv}
i\frac{\partial   L}{\partial   t^m_n} =
\frac{k_{n,m}}{\sqrt{r}}\left[(L^{n+\frac{m+1}{r}})_+,L\right],
\end{equation}
where the constants 
$$
k_{n,m}=\frac{(-1)^nr^{n+1}}{(m+1)(r+m+1)\ldots (nr+m+1)}
$$
have been inserted for convenience.

In terms of the unknown functions $u_k$ the equations~(\ref{eq:kdv}) take the 
form
$$
\frac{\partial u_k}{\partial  t^m_n} = S^m_{k,n},
$$
where $S^m_{k,n}$ are polynomials in $u_j$ and their derivatives with respect 
to $x$.

There are several simple special cases.
\begin{enumerate}

\item When $m+1$ is a multiple of $r$, the commutator in~(\ref{eq:kdv})
vanishes and the corresponding equation just means that the functions $u_k$ do 
not depend on $t^m_n$.

\item In the case $m=n=0$ the corresponding equation is 
$\frac{\partial u_k}{\partial  t^0_0}=
\frac{\partial u_k}{\partial  x}$ and therefore we can identify $t^0_0$ with $x$.

\item Finally, when $r=2$ the $\KdV_r$ hierarchy becomes the ordinary $\KdV$ hierarchy.
\end{enumerate}

\subsection{Potential}

In the notation of the previous subsection, we introduce the functions
\begin{equation}
  \label{eq:v}
  v_n = -\frac{r}{n+1} \mathrm{res}(L^{1/r})^{n+1},
\end{equation}
where the residue of a pseudodifferential operator is defined as the
coefficient of $D^{-1}$. The functions $v_k$ can be expressed in terms of
$u_j$ by a triangular system of differential polynomials. This means that
$u_j$ can be expressed in terms of $v_n$ in a similar way, and we can
consider $v_0,v_1,\ldots, v_{r-2}$ as a new system of coordinates in $\cd$.

We call a function $\Phit(\bt)$ in formal variables $t_n^m$, $n,m\ge 0$, a 
{\em   potential\/} of the $\KdV_r$ hierarchy if, first, $\Phit({\bf 0})=0$, 
second, the functions 
$$
v_m(\bt) = \frac{\partial \Phit(\bt)}{\partial  t^0_0 \partial t_0^m}
$$
satisfy the equations~(\ref{eq:kdv}) with $x=t_0^0$ and $u_j$ related with
$v_m$ via~(\ref{eq:v}) and, finally, $\Phit(\bt)$ satisfies the so-called 
{\em   string equation\/}
\begin{equation}
  \label{eq:string}
\frac{\partial \Phit(\bt)}{\partial  t^0_0 }  
= \frac{1}{2}\sum_{m,n=0}^{r-2} \eta_{mn}t_0^mt_0^n
+ \sum_{k=1}^\infty \sum_{m=0}^{r-2}t_{k+1}^m
\frac{\partial \Phit(\bt)}{\partial  t^n_m},
\end{equation}
where $\eta_{mn}=\delta_{m+n,r-2}$.

It can be shown that the potential $\Phit(\bt)$ is uniquely determined by these
conditions (cf.~\cite{W3}).

\subsection{Semiclassical approximation}

The hierarchy $\KdV_r$~(\ref{eq:kdv}) has a semiclassical (or dispersionless)
limit which is defined in terms of the formalism of Section~\ref{ssec:kdv} as follows.
For a differential operator $L\in \cd$ given
by~(\ref{eq:dops}) denote by $\widetilde{L}=p^r-\ds \sum_{m=0}^{r-2}u_m(x)p^m$
the polynomial in a formal variable $p$ obtained by replacing $D$ with $p$.
The commutator $[L,Q]$ of differential operators will be replaced
in~(\ref{eq:kdv}) by the Poisson bracket 
$$
\{\widetilde{L},\widetilde{Q}\} = 
\frac{\partial \widetilde{L}}{\partial   p}  \frac{\partial \widetilde{Q}}{\partial x} 
-\frac{\partial \widetilde{Q}}{\partial   p}  \frac{\partial \widetilde{L}}{\partial x} .
$$

{\em The semiclassical limit of the $\KdV_r$ hierarchy\/}, $\KdV^s_r$ is the system
of equations 
\begin{equation}
\label{eq:kdv0}
\frac{\partial \widetilde{L}}{\partial  t^m_n}
=\frac{k_{m,n}}{r} \left\{\widetilde{L^{n+\frac{m+1}{r}}},\widetilde{L}\right\}
\end{equation}
in unknown functions $u_0,\ldots,u_{r-2}$.

The corresponding potential function $\Phit_0(\bt)$ is defined as the unique
function satisfying the string equation~(\ref{eq:string}), the condition
$\Phit_0({\bf 0})=0$, and such that the functions $u_0,\ldots,u_{r-2}$ given
by~(\ref{eq:v}) and
$$
v_m(\bt) = \frac{\partial \Phit_0(\bt)}{\partial  t^0_0 \partial t_0^m}
$$
satisfy the equations of the hierarchy~(\ref{eq:kdv0}).

\section{The generalized Witten conjecture}

In \cite{W3}, Witten conjectured that the potential  $\tilde{\Phi}(\bt)$ of the
$\KdV$ hierarchy is equal to a generating function for the intersection
numbers of certain tautological classes on the moduli space of stable curves
$\M_{g,n}$. This conjecture arose from the study of physical realizations of
pure topological gravity and provided an unexpected link between integrable
systems and the geometry of these moduli spaces.  This conjecture was proved by
Kontsevich \cite{Ko}.

To be more precise, let $\M_{g,n}\,:=\,\{\,[\Sigma,
x_1,\ldots,x_n]\,\}$ be the moduli space of stable curves of genus $g$
with $n$ marked points. That is, $\Sigma$ is a conformal equivalence
class of genus $g$ Riemann surfaces, with (at worst) nodal
singularities and $x_1,\ldots, x_n$ are marked points, distinct from
each other and the nodes.  These data are subject to the stability
condition that they must have no infinitesimal automorphisms. This
moduli space is a ($3g-3+n$)-dimensional, compact, complex orbifold
and is a compactification of the moduli space of genus $g$ Riemann
surfaces with $n$ marked points.

This space $\mgnbar$ is equipped with  tautological holomorphic line bundles
$\cl_{(g,n),i}\,\to\,\M_{g,n}$ where $i\,=\,1,\ldots,n,$ such that the 
fiber of $\cl_{(g,n),i}$ at the point $[\Sigma; x_1,\ldots, x_n]$ is the
cotangent space $T_{x_i}^* \Sigma$. Define $\psi_i\,\in\,H^2(\M_{g,n})$ to be
$c_1(\cl_{(g,n),i})$ for all $i\,=\,1,\ldots, n$. Let
\[
\left< \tau_{a_1}\ldots \tau_{a_n}
\right>_g\,:=\,\int_{\M_{g,n}}\,\psi_1^{a_1}\ldots\psi_n^{a_n}
\]
for all nonnegative integers $a_1,\ldots,a_n\,\geq\,0$. Let
$\bt\,:=\,(t_0,t_1,\ldots)$ be formal variables and consider the 
function $\Phi_g(\bt)\,\in\,\nc[[\bt]]$ defined by
\begin{equation}\label{eq-phi}
\Phi_g(\bt)\,:=\,\left<\exp(\bt\cdot\btau)
\right>_g\,=\,\sum_{n\,\geq\,3}\,\frac{t_{a_1}\,\cdots\,t_{a_n}}{n!}\,\left<
\,\tau_{a_1}\cdots\tau_{a_n}\,\right>_g
\end{equation}
where $\bt\,\cdot\,\btau\,=\,\sum_{a=0}^\infty t_a \tau_a$, the
expression $\left<\cdots\right>_g$ is linear over $\nc[[\bt]]$, and
the exponential function is a formal power series. The summation
convention for repeated indices has also been used. Let
$\Phi(\bt)\,=\,\sum_{g\,\geq\,0}\,\Phi_g(\bt)$.

\begin{thm} (Kontsevich \cite{Ko})
  The function $\Phi(\bt)$ defined in (\ref{eq-phi}) is equal to the
  function $\Phit(\bt)$ for $\KdV_2$.
\end{thm}

Kontsevich's proof is fascinating, but it does not explain why such a
theorem \textsl{should} hold. Kontsevich realizes the function
$\tilde{\Phi}(\bt)$ as a large $N$ limit of an integral over the space
of Hermitian $N\,\times\,N$ matrices, a so-called Hermitian matrix
model. The terms in the resulting perturbative expansion are indexed
by ribbon graphs which, in turn, can be related to  a cell decomposition
(constructed using Strebel differentials) of a closely related moduli space.

A natural question to ask is whether there is a moduli space and
tautological classes on it, whose generating function assembles into
the potential $\Phit$ of the $\KdV_r$ hierarchy for general $r$.
Witten \cite{W,W2} stated a conjecture for this case as well. He
expected the relevant theory to be topological gravity coupled to a
topological matter sector corresponding to a class of topological minimal
models which he realized in terms of an
$\mathfrak{su}(2)_{r-2}/\mathfrak{u}(1)$ coset model. He suggested a
construction of the relevant moduli space and tautological classes.
However, the rigorous construction of the appropriate compact moduli
space had not yet been completed at the time that the conjecture was
formulated.

Recently, Jarvis \cite{J,J2} constructed the moduli spaces, called
\textsl{the moduli space of $r$-spin curves of genus $g$ with $n$
  marked points $\M_{g,n}^{1/r}$}. These spaces are the disjoint union
\[
\M_{g,n}^{1/r}\,=\,\sqcup_{\bm}\,\M_{g,n}^{1/r,\bm}
\]
where the union is over $n$-tuples of nonnegative integers
$\bm\,=\,(m_1,\ldots,m_n)$ where $0\,\leq\,m_i\,\leq\,r-1$ for all
$i\,=\,1,\ldots\,n$. (In fact, the moduli spaces are defined for general
integral $\bm$ but we will only consider those pieces relevant to the
generalized Witten conjecture.)

The spaces $\M_{g,n}^{1/r}$ are equipped with tautological 
classes $\psi_i\,\in\, H^2(\M_{g,n}^{1/r})$ where 
$i\,=\,1,\ldots,n$. Each $\psi_i$ is the Chern class of a 
tautological line bundle over $\M_{g,n}^{1/r}$. Furthermore, 
Witten \cite{W} sketched an analytic construction of a 
tautological class $\cv$ in $H^D(\M^{1/r}_{0,n})$ where 
$D\,=\,\frac{1}{r}(2-r+\sum_{i=1}^n m_i)$.  The correlators 
\[
\left<\,\tau_{a_1}(e_{m_1})\,\ldots\,\tau_{a_n}(e_{m_n})\,\right>_g\,:=\, r^{1-g}
\int_{\M_{g,n}^{1/r,\bm}}\, \cv\,\psi_1^{a_1}\,\ldots\,\psi_n^{a_n}
\]
assemble into a generating function $\Phi$, which is conjectured to
coincide with the potential $\Phit$ of $\KdV_r$.  The integral above
is to be understood in the orbifold sense.\footnote{The factor of $r$
  can be interpreted as coming from the string coupling constant.}

In genus zero a rigorous algebro-geometric construction of this 
class was given in \cite{JKV} (see also Section~\ref{cvirt}), and 
the following result holds. 

\begin{thm}(The Generalized Witten Conjecture in Genus Zero \cite{JKV})
\label{thm:gwconj}
Consider  the generating function
\[
\Phi_0(\bt)\,:=\,\left<\exp(\bt\cdot\btau)\,\right>_0\,:=\,\sum_{n\,\geq\,3}\,
\frac{t_{a_1}^{m_1}\cdots
t_{a_n}^{m_n}}{n!}\,\left<\,\tau_{a_1}(e_{m_1})\,\ldots\,\tau_{a_n}(e_{m_n})\, 
\right>_0 
\]
in $\nq[[\bt]]$ where $\bt$ is the set of variables $t_a^m$ for
$a\,\geq\,0$ and $m\,=\,0,\ldots,r-1$,
$\bt\,\cdot\,\btau\,:=\,\sum_{a,m} \tau_{(a,m)}\, t_a^m$, and
$a_i\,\geq\,0$ and $0\,\leq\,m_i\,\leq\,r-1$. The function
$\Phi_0(\bt)$ which is in fact independent of $t_a^{r-1}$ is 
equal to $\Phit_0(\bt)$, the (semiclassical) potential of the 
$\KdV_r^s$ hierarchy. 
\end{thm}

The precise algebro-geometric formulation of this conjecture is, 
as yet, incomplete in higher genera since the class $\cv$ has not 
yet been constructed algebro-geometrically 

More detailed descriptions of the above objects will be given in
subsequent sections.

\section{The moduli space of curves with multiple spin structures}

In this section we review some definitions and facts from \cite{JKV}
about spin curves, and then generalize these to define curves with
multiple spin structures.

\subsection{Definitions}

\begin{df}
Let $(X, p_1, \dots, p_n)$ be a nodal, $n$-pointed algebraic 
curve, and let $\ck$ be a rank-one, torsion-free sheaf on $X$. A 
\emph{$d$-th root of $\ck$ of type $\mathbf{m}=(m_1, \ldots, 
m_n)$} is a pair $(\ce, b)$ of a rank-one,  torsion-free sheaf 
$\ce$, and an $\co_X$-module homomorphism 
$$ 
  b: \ce^{\tensor d}  \rTo \ck \otimes \co_X(-\sum m_ip_i)
$$ 
with the following properties:
\begin{itemize}
\item $d \cdot \deg \ce = \deg \ck-\sum m_i$
\item $b$ is an isomorphism on the locus of $X$ where $\ce$ is 
locally free
\item for every point $p \in X$ where $\ce$ is not free, the 
length of the cokernel of $b$ at $p$ is $d-1$.
\end{itemize}
\end{df}

For any $d$-th root $(\ce,b)$ of type $\mathbf{m}$, and for any 
$\mathbf{m}'$ congruent to $\mathbf{m} \mod d$, we can construct 
a unique $d$-th root $(\ce',b')$ of type $\mathbf{m}'$ simply by 
taking $\ce'=\ce \otimes \co(1/d \sum (m_i-m_i')p_i)$. 
Consequently, the moduli of curves with $d$-th roots of a bundle 
$\ck$ of type $\mathbf{m}$ is canonically isomorphic to the 
moduli of curves with $d$-th roots of type $\mathbf{m}'$.  
Therefore, unless otherwise stated, we will always assume the 
type $\mathbf{m}$ of a $d$-th root has the property that $0 \leq 
m_i<d$ for all $i$.

Unfortunately, the moduli space of curves with $d$-th roots of a 
fixed sheaf $\ck$ is not smooth when $d$ is not prime, and so we 
must consider not just roots of a bundle, but rather coherent 
nets of roots \cite{J}.  This additional structure suffices to 
make  the moduli  space of curves with a coherent net of roots 
smooth.

\begin{df}
Let $\ck$ be a rank-one, torsion-free sheaf on a nodal 
$n$-pointed curve $(X, p_1, \ldots, p_n)$. A \emph{coherent net 
of $r$-th roots of $\ck$ of type $\mathbf{m}=(m_1, 
\ldots, m_n)$} consists of the  following data: 
\begin{itemize}
\item For every divisor $d$ of $r$, a rank-one 
torsion-free sheaf $\ce_d$ on $X$; 

\item For every pair of divisors $d',d$ of $r$,  such that $d'$ divides $d$,  
an $\co_X$-module  homomorphism 
$$
c_{d,d'} : \ce^{\tensor d/d'}_{d} \rTo  
\ce_{d'}. 
$$
\end{itemize}

These data are subject to the following restrictions.
\begin{enumerate}
\item $\ce_1=\ck$ and $c_{1,1}=\mathbf{1}$
\item For each divisor $d$ of $r$ and each divisor $d'$ of $d$, 
let $\mathbf{m}''=(m''_1, \ldots, m_n'')$ be such that $m''_i$ is 
the unique non-negative integer less than $d/d'$, and congruent 
to $m_i \mod d$. Then the  homomorphism $c_{d,d'}$  makes $(\ce_d, 
c_{d,d'})$ into a $d/d'$ root of $\ce_{d'}$ of type $\mathbf{m}''$. 
\item The homomorphisms $\{c_{d,d'}\}$ are compatible.  
That is, the diagram 
$$
\begin{diagram}
(\ce^{\tensor d/d'}_{d})^{\tensor d'/d''} & \rTo^{(c_{d,d'})^{\tensor d'/d''}} 
& \ 
\ce^{\tensor d'/d''}_{d'}\\ 
&  \rdTo^{c_{d,d''}}  & \dTo_{c_{d',d''}} \\
& &  \ce_{d''}\\
\end{diagram}
$$
commutes for every $d''|d'|d|r$.  
\end{enumerate} 
\end{df}

If $r$ is prime, then a coherent net of $r$-th roots is simply an 
$r$-th root of $\ck$.  Moreover, if $\ce_d$ is locally free of 
type $\bm$, then 
for $d'|d$ if $\bm \equiv \bm \pmod d'$ and $0 \leq m'_i <d'$, 
then the sheaf $\ce_{d'}$ is isomorphic to $\ce_{d}^{\tensor 
d/d'}\otimes \co(\frac{1}{d'} \sum (m_i -m'_i)p_i)$ and 
$c_{d,d'}$ is uniquely determined up to an automorphism of 
$\ce_{d'}$. 

\begin{df}
An \emph{$n$-pointed $r$-spin curve of type $\bm = (m_1, \ldots, 
m_n)$} is an $n$-pointed, nodal curve $(X, p_1, \ldots, p_n)$ 
with a coherent net of $r$-th roots of $\omega_X$ of  type 
$\mathbf{m}$, where $\omega_X$ is the canonical (dualizing)  
sheaf of $X$.  An $r$-spin curve is called \emph{smooth},  if $X$ 
is smooth, and it is called \emph{stable}, if $X$ is stable.  
\end{df}

\begin{ex}
Smooth 2-spin curves of type $\bm =\mathbf{0}$
correspond to classical spin curves (a curve and a 
theta-characteristic) \emph{together with an explicit isomorphism 
$\ce^{\tensor 2}_{2}  \irightarrow \omega$.} 
\end{ex}

\begin{df}
An isomorphism of $r$-spin curves 
$$
(X, p_1, \ldots, p_n, \{\ce_d, c_{d, d'}\}) 
\irightarrow (X', p'_1, \ldots, p'_n, \{\ce'_d, c'_{d,  d'}\})
$$ 
of the same type $\bm$ is an isomorphism of pointed curves  
$$ 
  \tau : (X, p_1,  \ldots, p_n)  \irightarrow (X', p'_1, \ldots, p'_n)
$$ 
and a family of 
isomorphisms $\{\beta_d : \tau^* \ce'_d  \irightarrow \ce_d\},$ 
with $\beta_1$  being the canonical isomorphism $\tau^* 
\omega_{X'}(-\sum_i m_i  {p'}_i)  \irightarrow 
\omega_X(-\sum m_ip_i),$ and such that the $\beta_d$ are compatible with 
all the maps $c_{d,d'}$ and $\tau^*c'_{d,d'}$.
\end{df}

Every $r$-spin structure on a smooth curve $X$ is determined up 
to isomorphism by a choice of a line bundle $\ce_r$, such that 
$\ce^{\tensor r}_r \cong \omega_X (-\sum m_i p_i)$.  In 
particular, if $X$ has no automorphisms, then the set of 
isomorphism classes of $r$-spin structures of a fixed type 
$\mathbf{m}$ is a principal $\jac_rX$-bundle, where $\jac_rX$ is 
the group of $r$-torsion points of the Jacobian of $X$, thus 
there are $r^{2g}$ such isomorphism classes.

\begin{ex}
  If $g=1$ and $\bm=\mathbf{0}$, then $\omega_X \cong \co_X$ and a
  smooth $r$-spin curve is just an $r$-torsion point of the elliptic
  curve $X$, again with an explicit isomorphism of $\ce^{\tensor
    r}_{r} \irightarrow \co_X$.  In particular, the stack of stable
  $r$-spin curves of genus one and type $\mathbf{0}$ forms a gerbe
  over the disjoint union of modular curves $\coprod_{d|r} X_1(d)$.
\end{ex}
\begin{df}

The moduli space of stable $n$-pointed, $r$-spin curves of genus 
$g$ and type $\mathbf{m}$ is denoted $\mgnrmbar$.  The disjoint 
union $\displaystyle\coprod_{\substack{\mathbf{m} \\ 0 \leq m_i 
<r}}  \mgnrmbar$ is denoted $\mgnrbar$.    
\end{df}

\begin{rem} \label{rem:restrict}
  As mentioned above, no geometric information is lost by restricting
  $\bm$ to the range $0\le m_i \le r-1$, since when $\bm \equiv \bm'
  \mod r$ every $r$-spin curve of type $\bm$ naturally gives an
  $r$-spin curve of type $\bm'$ simply by
$$
\ce_d \mapsto \ce_d \tensor \co (\sum 
\frac{m_i-{m'}_i}{d} p_i),
$$ 
and thus $\mgnrmbar$ is canonically isomorphic 
to  $\overline{\scheme{M}}^{1/r,\mathbf{m}'}_{g,n}$. 
\end{rem}

More generally, if $\mathbf{r} = (r^1, \dots, r^k)$ is a 
$k$-tuple of integers $r^i \geq 2$, and if  $\bbm 
:= (\bm^1,\ldots,\bm^k)$, $\bm^j := (m_1^j,\ldots,m_n^j)$, for 
nonnegative integers $m_i^j$ we define $\mathbf{r}$-spin curves 
of type $\bbm$ as follows: 
\begin{df}
An $n$-pointed $\mathbf{r}$-spin curve of type $\bbm$ is an 
$n$-pointed, nodal curve $(X, p_1, \dots, p_n)$ with a $k$-tuple 
of $\mathbf{r}$-spin structures of type $\bbm$, that is, a 
$k$-tuple whose $i$-th term is a coherent net 
$\mathcal{N}_i:=\{\ce_{d^i}, c_{d^i, d^{i'}}\}$ of $r^i$-th roots 
of $\omega_X$ of type $\bm^i =(m^i_1, \dots, m^i_n)$.  Here, of 
course, the $d^i$ run over all positive integers dividing $r^i$, 
and the $d^{i'}$ run over all positive divisors of $d^i$. 
\end{df}

\begin{df}
  We denote by $\M^{\mathbf{1/r}, \mathbf{m}}_{g,n}$ the moduli of
  genus $g$, $n$-pointed, stable $\mathbf{r}$-spin curves of type
  $\bbm$ and we denote by $\M^{\mathbf{1/r}}_{g,n}$ the disjoint union
  $$
  \M^{\mathbf{1/r}}_{g,n}:= \coprod_{\bbm=(\bm^1, \dots,\bm^k)
    \atop {\bm^i =(m^i_1, \dots, m^i_n) \atop 0 \leq m^i_j<r^i}}
  \M^{\mathbf{1/r},\bbm}_{g,n}$$
\end{df}

\subsection{Basic Properties of the Moduli Space $\M^{\mathbf{1/r},
    \bbm}_{g,n}$}
These spaces are endowed with a number of canonical morphisms,
including the projections $$q_i:\M^{\mathbf{1/r}, \bbm}_{g,n} \rTo
\M^{1/r^i,\bm^i}_{g,n}$$
and $$p_i:\M^{1/r^i,\bm^1}_{g,n} \rTo
\M_{g,n}.$$
It is clear that
\begin{equation}\label{fiber-prod}
\M_{g,n}^{\bo/\br,\bbm} = \prod_{i = 1}^k \M_{g,n}^{1/r^i,\bm^i}
\end{equation}
where $\prod_i$ denotes the fibered product over $\M_{g,n}$  with 
respect to the maps $p_i$.

In \cite{J} it is shown that $\mgnrmbar$ is a smooth Deligne-Mumford
stack (orbifold), finite over $\mgnbar$, with a projective coarse
moduli space. Consequently, $\M^{\mathbf{1/r}, \bbm}_{g,n}$ has the
same properties.

There is one other canonical morphism associated to these spaces; 
namely, when $d$ divides $r$, the morphism 
$$ [r/d]:\mgnrmbar \rTo \mbar^{1/d,\mathbf{m}''}_{g,n}
\mathrm{\ and\ } 
[r/d]:\mgnrbar \rTo \mbar^{1/d}_{g,n},
$$
which forgets all of the roots and homomorphisms in the net of 
$r$-th roots except those associated to divisors of $d$.  Here, as above, 
$\mathbf{m}''$ is   congruent to $\mathbf{m} \mod d$ and has $0 
\leq m''_i <d$.  In the case that the underlying curve is smooth, 
this is equivalent to replacing the line bundle $\ce_r$ by its 
$r/d$-th tensor power (and tensoring with $\co(1/d \sum(m_i-m''_i)p_i )$). 

\begin{rem}
  It is important to keep in mind that the space $\M^{1/rs}_{g,n}$ is
  not isomorphic to $\M^{(1/r,1/s)}_{g,n} =\M^{1/r}_{g,n}
  \cross_{\M_{g,n}} \M^{1/s}_{g,n}$, even when $r$ and $s$ are
  relatively prime.  However, on the open locus of smooth spin curves
  we do have $\mathcal{M}^{1/rs}_{g,n}\cong
  \mathcal{M}^{(1/r,1/s)}_{g,n},$ whenever $r$ and $s$ are relatively
  prime \cite{J}.
\end{rem}
Throughout this paper we will denote the universal curve  by 
$\pi:\mathcal{C}^{\mathbf{1/r}}_{g,n} \rTo 
\M^{\mathbf{1/r}}_{g,n}$.  As in the case of the moduli space of 
stable curves, the universal curve possesses canonical sections 
$\sigma_i:\M^{\mathbf{1/r}}_{g,n} 
\rTo \mathcal{C}^{\mathbf{1/r}}_{g,n}$ for $i = 1 \ldots n$. 

\subsection{Cohomology Classes}\label{cvirt}
As in the case of stable maps and stable curves we have some 
tautological cohomology classes.  Most important for our purposes  
are the classes $$\psi_i 
:=c_1(\cl_{(g,n),i}) = c_1(\sigma^*_i(\omega_{\pi})),$$ the first Chern class 
of the restriction of the canonical (relative dualizing) sheaf 
$\omega_{\pi}$ to the image of the $i$-th section.  Of course, 
$\psi_i$ on $\M^{\mathbf{1/r},\bm}_{g,n}$ is just the pullback 
$p^*\psi_i$ of the usual class $\psi_i$ on $\M_{g,n}$ because 
$\cl_{(g,n),i} = \sigma_i^*(\omega_{\pi}) = p^* 
\sigma_i^*(\omega_{(g,n)})$, where $\omega_{(g,n)}$ is the relative 
dualizing sheaf of the universal curve over $\M_{g,n}$.

In addition to the tautological classes, we will also need a 
class $c^{1/r}$ which plays the role of Gromov-Witten classes 
$(ev^*_1(\gamma_1)\cup \dots \cup ev^*_n (\gamma_n)) \cap 
([\M_{g,n}(V,\beta)]^{virt}))$ in the case of stable maps.

The class $c^{1/r}$ will be called the {\em virtual class}, and in the
case of $g=0$ it has a simple description as the Euler class of the
pushforward of the $r$-th root bundle from the universal curve to
$\mgnrbar$.

\begin{df}
If $g=0$, let $\ce^i_{r^i}$ be the $r^i$-th root of $\omega_{\pi}$ 
associated to the $i$-th net $\mathcal{N}_i$ of the universal 
$\mathbf{r}$-spin structure $(\mathcal{N}_1, \dots, 
\mathcal{N}_k)$ on the universal curve $\pi: \mathcal{C}^{\mathbf{1/r}}_{0,n}
\rTo \M^{\mathbf{1/r}}_{0,n}$.  Let $c^{\mathbf{1/r}}_{0,n}$ be the Euler
class (top Chern class) of the pushforward of $\ce^1_{r^1}\oplus 
\dots \oplus \ce^k_{r^k}$ 
$$c^{\mathbf{1/r}}_{0,n} = e(-R\pi_*(\ce^1_{r^1}
\oplus \dots \oplus \ce^k_{r_k})).$$ 
\end{df}

Here $R \pi_*\cf$ is the $K$-theoretic pushforward $\pi_* 
\cf-R^1\pi_*\cf$.  This class is well-defined because in genus 
zero $\ce^i_{r^i}$ is {\em convex} \cite{JKV}, that is $\pi_* 
\ce^i_{r^i}=0$, and so $R^1 \pi_* \ce^i_{r^i}$ is a vector bundle 
of rank $D_i =\chi(\ce^i_{r^i})= ((2-r_i)+\sum_j m^i_j)/r_i$.  
Thus $R \pi_*(\ce^1_{r^1} \oplus \dots \oplus \ce^k_{r^k})$ is 
also a vector bundle and has a top Chern class. 
 
Let us consider the case where $k=1$ and $r \geq 2$. The class $\cv$
on $\M_{0,n}^{1/r}$ defined above is precisely the class described in
Theorem~\ref{thm:gwconj} (the Generalized Witten Conjecture).

\begin{lm} \label{pb-prod} 
  When the moduli space $\M_{g,n}^{\bo/\br}$ of curves with multiple
  spin structures is viewed as a fibered product, as in
  (\ref{fiber-prod}), then the virtual class $c^{\bo/\br}$ is the
  product of the pullbacks of the virtual classes from each of the
  factors.
  $$c^{\mathbf{1/r}} = q_1^* c^{1/r^1} \cup q_2^* c^{1/r^2} \cup \dots
  q_n^* c^{1/r^n}$$

\end{lm}

\begin{proof}
Since the maps $p_i:\M^{1/r^i}_{g,n} \rTo \M_{g,n}$ are flat, 
since $\M^{\mathbf{1/r}}_{g,n}$ is a fibered product, and since
$\ce^i_{r^i}$ is the pullback of the universal $r^i$-th root on 
$\mathcal{C}^{1/r^i}_{0,n}$ via $q_i$, we have 
\begin{eqnarray*}
c^{\mathbf{1/r}} & = & e(-(R\pi_*(\ce^1_{r^1} \oplus \dots \oplus 
\ce^k_{r^k}))) \\
&=& e(-(q^*_1 R\pi_* \ce^1_{r^1} \oplus \dots \oplus q^*_k R 
\pi_* 
\ce^k_{r^k}))\\
&=& q_1^* c^{1/r^1} \cup q_2^* c^{1/r^2} \cup \dots q_n^* c^{1/r^n}
\end{eqnarray*} 
\end{proof}

\begin{rem}
In \cite{JKV} it is shown that in genus zero the virtual class 
$\cv$ has a number of properties resembling those of the 
Gromov-Witten classes.  One additional property that is 
unique to this construction, and which plays an important role, is 
the fact that $\cv$ vanishes on $\M_{g,n}^{\bo/\br,\bm}$ 
whenever any one of the $m_j^i$ is equal to $r^i-1$.  Also it
vanishes in some other special cases.
\end{rem}

\section{Cohomological Field Theories and Frobenius Manifolds}
\subsection{Operads and Cohomological Field Theories}

The objects defined in the previous section allow a construction of a $g=0$
cohomological field theory (\cft) whose potential function coincides
with the potential function of $\KdV_r$ and gives rise to an
associated $(r-1)$ dimensional Frobenius manifold. 
It is precisely this property, together with the topological
recursion relations (TRR) which was used to prove the generalized
Witten conjecture in genus zero. It is also the framework in which the
notion of tensor product 
most readily appears.  These constructions fit into the framework of stable
graphs and operads which we will now describe.

Consider a connected graph $\Gamma$ with $n$ tails (or legs) indexed
by numbers $1,\ldots, n$, such that each vertex $v$ is
labeled by a nonnegative integer $g(v)$ (called the \textsl{genus of 
  $v$}).
 The \textsl{genus}  $g(\Gamma)$ of $\Gamma$ is defined as
$g(\Gamma)=\sum_{v\in V(\Gamma)}g(v) + b_1(\Gamma)$, where 
$V(\Gamma)$ is the set of vertices of $\Gamma$ and $b_1(\Gamma)$ 
is the first Betti number (number of loops) of $\Gamma$. Such a 
graph $\Gamma$ is called a \textsl{stable graph} of genus $g$ 
with $n$ tails if for every $v\in V(\Gamma)$ with $g(v)=0$ its 
valence (denoted by $n(v)$) satisfies $n(v)\,\geq\,3$ and if 
$g(v)=1$ then $n(v)\,\geq\,1$. In particular,  genus zero stable 
graphs are trees whose vertices are all labeled by $0$.

Each  stable curve $[\Sigma,x_1,\ldots,x_n] \in \M_{g,n}$
is associated to a stable graph of
genus $g$ with $n$ tails in the following way. 

To each irreducible component, associate a vertex, and draw as
many half-edges emerging from this vertex as there are \textsl{special
  points} (either marked points or nodes)  on that component.
 Label vertex $v$ by an integer $g(v)$ equal to the genus of that
 irreducible  component. Connect any two half 
edges if they are associated to the same node. Finally, label the 
label the remaining half-edges (tails) by an 
integer from $1,\ldots,n$ corresponding to the associated marked point.

Let $\Gamma$ be a stable graph of genus $g$ with $n$ tails. The
stratification of $\M_{g,n}$ yields \textsl{gluing maps} 
\[
\rho_\Gamma\,:\,\prod_{v\,\in\,V(\Gamma)}\,
\M_{g(v),n(v)}\,\to\,\M_{g,n}. 
\] 
This map is obtained by composing the maps 
$\chi\,:\,\prod_{v\,\in\,V(\Gamma)}\,\M_{g(v),n(v)}\,\to\,\M_\Gamma$ 
and $i\,:\,\M_\Gamma\,\to\,\M_{g,n}$.  The space $\M_\Gamma$ is 
the closure of the locus of points in $\M_{g,n}$ whose associated 
stable graph is $\Gamma$. The map $\chi$ is the quotient map 
corresponding to the action of $\Aut(\Gamma)$, the automorphism 
group of $\Gamma$, and $i$ is the inclusion map. The symmetric 
group $S_{n(v)}$ acts on the space $\M_{g(v),n(v)}$ by permuting 
the marked points. Similarly, $S_n$ acts on  $\M_{g,n}$. Clearly, 
the map $\rho_\Gamma$ is equivariant under these actions. 

The collection $\M \, :=\, \{\,\M_{g,n}\,\}$, together with the 
symmetric group actions and the maps $\rho_\Gamma$, is the model 
example of a \textsl{modular operad} \cite{GK}. Restricting to 
genus zero stable graphs, one obtains a \textsl{cyclic operad}, a 
refinement of the notion of an operad. The modular operad 
structure on $\M$ induces a modular operad structure on the 
homology $H_\bullet(\M)\,:=\,\{ H_\bullet(\M_{g,n}) \}$. From 
this perspective, a cohomological field theory (\cft) is simply a 
vector space $\ch$ with a metric $\eta$,  which is an algebra 
over $H_\bullet(\M)$.  The algebra maps 
$H_\bullet(\M_{g,n})\,\to\, T^n \ch^*$ (where $T^n 
\ch^*$ denotes the $n$-fold tensor product of $\ch^*$)  are the 
correlators of the theory. However, in the context of algebraic 
geometry, it turns out to be more useful to use the dual 
definition of a cohomological field theory, which is given in 
terms of cohomology. 

\begin{df}
A \emph{cohomological field theory (\cft) of rank $d$} (denoted 
by $(\ch,\eta,\Lambda)$ or just $(\ch,\eta)$) is a 
$d$-dimensional vector space $\ch$ with a metric $\eta$ and a 
collection of forms $\Lambda\,:=\, 
\{\,\Lambda_{g,n}\,\}$ where
\begin{equation}\label{eq:lambda}
\Lambda_{g,n} \in \,H^\bullet(\mgnbar)\,\otimes T^n \ch^*
= Hom(T^n \ch,H^\bullet(\mgnbar))
\end{equation} 
are defined for stable pairs $(g,n)$ and satisfy the following 
axioms {\bf C1--C3}.  Let $\{e_0,\ldots,e_{d-1}\}$ be a fixed 
basis of $\ch$, and let $\eta^{\mu\nu}$ be the inverse of the 
matrix of the metric  $\eta$ in this basis.  We use the summation 
convention in the following. 

\begin{enumerate}
\item[\bf C1.] The form 
$\Lambda_{g,n}$ is invariant under the action of the symmetric
  group $S_n$.
\item[\bf C2.] 
Let 
\begin{equation}
\label{eq:gluetree}
\rho_{\Gamma_\mathrm{tree}}:
\M_{k,j+1} \times \M_{g-k,n-j+1} \rTo \M_{g,n}
\end{equation} 
be the gluing map corresponding to the stable graph
$$
\Gamma_\mathrm{tree}=\tcongi
$$
then the forms $\Lambda_{g,n}$ satisfy the composition property 
\begin{eqnarray}
\label{eq:cfttree}
\rho_{\Gamma_\mathrm{tree}}^* \Lambda_{g,n}(\gamma_1,\gamma_2\,\ldots,\gamma_n)&=&
\\ 
\Lambda_{k,j+1}(\gamma_{i_1},\ldots,\gamma_{i_j},e_\mu)
\eta^{\mu\nu} &\otimes&
\Lambda_{g-k,n-j+1}(e_\nu\,\gamma_{i_{j+1}},\ldots,\gamma_{i_n})
\nonumber
\end{eqnarray} 
for all $\gamma_i \in \ch$.

\item[\bf C3.] 
Let 
\begin{equation}
\label{eq:glueloop}
\rho_{\Gamma_\mathrm{loop}}:
\M_{g-1,n+2} \rTo \M_{g,n} 
\end{equation}
be the gluing map corresponding to the stable graph
\begin{equation}
  \label{eq:gtree}
\Gamma_\mathrm{loop}=\ocongi
\end{equation}
then
\begin{equation}
\label{eq:cftloop}
\rho_{\Gamma_\mathrm{loop}}^*\,\Lambda_{g,n}(\gamma_1,\gamma_2,\ldots,\gamma_n)\,=\, 
\Lambda_{g-1,n+2}\,(\gamma_1,\gamma_2,\ldots,\gamma_n, e_\mu,
e_\nu)\,\eta^{\mu\nu}. 
\end{equation}
\medskip

The pair $(\ch,\eta)$ is called the \emph{state space} of the \cft. The maps
$H_\bullet(\mgnbar)\, \to\, T^n \ch^*$ given by
$[c]\,\mapsto\,\int_{[c]}\Lambda_{g,n}$ are called the \emph{correlators of
  the \cft}. 
\medskip

An element $e_0\in \ch$ is called \emph{a flat identity}
of the \cft, if in addition, the following equations hold.

\item[\bf C4a.] 
For all $\gamma_i$ in $\ch$ we have
\begin{equation}
  \label{eq:identity}
\Lambda_{g,n+1}(\gamma_1,\ldots, \gamma_n, e_0)\, = 
\,\pi^*\Lambda_{g,n}(\gamma_1,\ldots, \gamma_n),
\end{equation}
where $\pi:\M_{g,n+1} \to \M_{g,n}$ is the universal curve on $\mgnbar$ and
\item[\bf C4b.]
\begin{equation}
\label{eq:identity2}
\int_{\M_{0,3}}\,\Lambda_{0,3}(\gamma_1,\gamma_2,e_0) = \eta(\gamma_1,\gamma_2)
\end{equation}
\end{enumerate}
A \cft\ with flat identity is denoted by $(\ch,\eta,\Lambda,e_0).$ 
\medskip

A \emph{genus ${g}$ \cft} on the state space $(\ch,\eta)$ is the 
collection of forms $\{\,\Lambda_{\tilde{g},n}\,\}_{\tilde{g}\leq 
{g}}$ that satisfy only those of the 
equations~(\ref{eq:cfttree}), (\ref{eq:cftloop}), 
(\ref{eq:identity}), and (\ref{eq:identity2}), where 
$\tilde{g}\le {g}$. 
\end{df}

\begin{rem}
The state space $\ch$ of a \cft\ can, in general, be 
$\mathbb{Z}_2$-graded, but for simplicity we will assume that 
$\ch$ contains only even elements as is the case for $\KdV_r$. 
\end{rem}

Let $\Gamma$ be a genus $g$ stable graph with $n$ tails, then, 
since the map $\rho_\Gamma$ can be constructed from gluing 
morphisms~(\ref{eq:gluetree}) and (\ref{eq:glueloop}), the forms 
$\Lambda_{g,n}$ satisfy the restriction property 
\begin{equation} \label{restriction}
\rho_\Gamma^*\Lambda_{g,n}\,=\, \eta_\Gamma^{-1}(\,\bigotimes_{v\in V(\Gamma)}
\Lambda_{g(v),n(v)}\,)
\end{equation}
where \[ \eta^{-1}_\Gamma\,:\,\bigotimes_{v\in V(\Gamma)} 
T^{n(v)}\,\ch^* 
\,\to\, T^n\,\ch^* \] 
contracts the factors $T^n\ch^*$ by means of the inverse of the metric $\eta$
and successive application of equations~(\ref{eq:cfttree})
and~(\ref{eq:cftloop}) . 

\begin{df}
The {\em potential function} of the \cft\
$(\ch,\eta,\Lambda)$ is a formal series $\Phi\,\in\,\nc[[\ch]]$ given by
\begin{equation}
   \label{eq:small}
\Phi(\bx)\, := \,\sum_{g\,=\,0}^\infty\,\Phi_g(\bx),
\end{equation}
where 
\[
\Phi_g(\bx)\,:=\,\sum_n\,\frac{1}{n!}\,\int_{\mgnbar}\,
\left<\,\Lambda_{g,n}\,,\,\bx^{\otimes n}\, \right>.
\]
Here $\left<\,\cdots\,\right>$ denotes evaluation, the sum 
over $n$ is understood to be over the stable range, and 
$\bx\,=\,\sum_{\alpha}\,x^\alpha\,e_\alpha$, where 
$\{\,e_\alpha\,\}$ is a basis of $\ch$.\footnote{The \textsl{string coupling
constant} factor of $\lambda^{2 g-2}$ in front of $\Phi_g(\bx)$ has been
suppressed to avoid notational clutter.}
\end{df}

All of the information of a genus zero \cft\  is encoded in its 
potential. 

\begin{thm} \cite{KM1,Ma} 
\label{thm:wdvvg}
An element $\Phi_0$ in $\nc[[\ch]]$ is the potential of a rank 
$d$, genus zero \cft\ $(\ch,\eta)$ if and only if it contains 
only terms which are of cubic and higher order in the coordinates 
$x^0,\ldots,x^{d-1}$ (corresponding to a basis 
$\{\,e_0,\,\ldots\,e_{d-1}\,\}$ of $\ch$) and it satisfies the 
{\em associativity, or WDVV} (Witten-Dijkgraaf-Verlinde$^2$) 
equation 
\[
\label{eq:WDVV}
\partial_{a} \partial_{b} \partial_{e}\Phi_0\, \eta^{ef}\, \partial_{f}
\partial_{c} \partial_{d}\Phi_0\, = \, \partial_{b} \partial_{c}
\partial_{e} \Phi_0\, \eta^{ef} \,\partial_{f} \partial_{a}
\partial_{d}\Phi_0,
\]
where $\eta^{ab}$  is the inverse of the matrix of $\eta$ in the 
basis $\{e_a\}$,  $\partial_a$ is derivative with respect  to 
$x^a$, and the summation convention has been used.

Conversely,  a genus zero \cft\ structure on $(\ch,\eta)$ is 
uniquely determined by its potential $\Phi_0$, which must satisfy 
the WDVV equation. 
\end{thm} 

A formal power series $\Phi$ which satisfies the WDVV equations  
determines a \emph{formal Frobenius manifold}, so the theorem 
shows that a genus zero \cft\ with flat identity is equivalent to 
endowing the state space $(\ch,\eta)$ with the structure of a 
formal Frobenius manifold \cite{Du,Hi,Ma}. The theorem follows 
from the  work of Keel \cite{Ke} who proved that 
$H^\bullet(\M_{0,n})$ is generated by boundary classes, and that 
all relations between boundary divisors arise from lifting the 
basic codimension one relation on $\M_{0,4}$. 

\subsection{Gromov-Witten invariants}
The best known examples of  \cfts\  come from the Gromov-Witten 
classes of a smooth, projective variety $V$, where the state 
space $\ch$ is $H^\bullet(V)$, and $\eta$ is the Poincar\'e 
pairing. Let $\M_{g,n}(V)$ be the moduli space of stable, genus 
$g$ maps into $V$ with $n$ marked points, a compactification of 
the moduli space of holomorphic maps $f:\Sigma \to V$ from a 
genus $g$ Riemann surface $\Sigma$ with $n$ marked points 
$(x_1,\dots,x_n)$ into $V$. There are canonical evaluation maps 
$\ev_i:\M_{g,n}(V)\,\to\,V$ corresponding to evaluating $f$ at 
$x_i$.  There is also the stabilization morphism 
$\mathrm{st}:\M_{g,n}(V)\,\to\,\M_{g,n}$ which ``forgets'' $f$. 
Finally, the space $\M_{g,n}(V)$ has a virtual fundamental class 
$[\M_{g,n}(V)]^\mathrm{virt}$ in the Chow group of $\M_{g,n}(V)$. 
In this situation, $(H^\bullet(V),\eta,\Lambda,\mathbf{1})$ is a 
\cft\  with flat identity, where $\bo$ is the usual unit 
element and 
\[
\Lambda_{g,n}(\gamma_1,\ldots,\gamma_n)\,:=\,
\mathrm{st}_*(\ev_1^*\gamma_1\cup\ldots\cup\ev_n^*\gamma_n\cap 
[\M_{g,n}(V)]^\mathrm{virt})
\]
for $\gamma_i\,\in\,H^\bullet(V)$. In the case of $g=0$, the 
correlators endow $H^\bullet(V)$ with a multiplication , which is 
a deformation  of the cup product known as the quantum cohomology 
of $V$. 

The \emph{large phase space potential} $\Phi(\bt)$ is the sum of 
formal power series $\Phi_g(\bt)$ for $g\,\geq\,0$ where 
\[
\Phi(\bt)_g\,:=\,\left<\,\exp(\bt\cdot\btau)\,\right>_g\,=\,\sum\,
\frac{t_{a_1}^{m_1}\cdots t_{a_n}^{m_n}}{n!}
\left< \tau_{a_1}(e_{m_1})\cdots\tau_{a_n}(e_{m_n})\right>_g
\]
for a given basis $\{\,e_m\,\}$ of $H^\bullet(V)$ and
\[
\left< \tau_{a_1}(e_{m_1})\cdots\tau_{a_n}(e_{m_n})\right>_g\,=\,
\int_{[\M_{g,n}(V)]^\mathrm{virt}}\,\ev_1^* e_{m_1}\psi_1^{a_1}\ldots\ev_n^*
e_{m_n} \psi_n^{a_n}.
\]
Setting $t_a^m$ to zero and $x^m := t_0^m$, one recovers the 
potential $\Phi(\bx)$ of the \cft\ associated to $\Lambda$ which 
is the generating function for the Gromov-Witten invariants of 
$V$.\footnote{To avoid notational clutter, factor of $q^\beta$ 
where $\beta\in H_2(V,\nz)$ in the Novikov ring has been 
suppressed.} 

If $V$ is a point then $\M_{g,n}(V)$ reduces to $\M_{g,n}$ and 
one has Kontsevich's theorem, which identifies $\Phi$ with the 
potential of the $\KdV_2$ hierarchy.   For general $V$ a similar 
relation is expected, but the corresponding integrable system has  
only been described in $g=0$ \cite{Du}, and in some cases, in 
$g=1$ \cite{DuZh}. However, there is a related conjecture by 
Eguchi, Hori, Xiong, and S.~Katz  \cite{EHX}, which states that 
the exponential of $\Phi(\bt)$ satisfies a highest weight 
condition for an action of the Virasoro algebra. The evidence for 
this conjecture is growing \cite{FaPa, FaPa2, GePa}. 

\subsection{The $KdV_r$ Frobenius manifold}
Consider an integer $r\,\geq\,2$. Let $(\chr,\eta)$ be an 
$(r-1)$-dimensional vector space $\chr$ with basis 
$\{\,e_0,\ldots\,e_{r-2}\,\}$ and a metric 
$\eta(e_{m'},e_{m''})\,=\,\delta_{m'+m'',r-2}$ such that 
$m',m''\,=\,0,\ldots,r-2$. 

\begin{thm}\cite{JKV}
For each integer $r\,\geq\,2$, the pair $(\chr,\eta)$ is a 
Frobenius manifold (and thus defines a $g=0$ \cft) with flat 
identity $e_0$, Euler vector field 
\[
E\,:=\,\sum_{m\,=\,0}^{r-2}\,(\frac{m}{r}-1)\,x^m\,\frac{\partial}{\partial x^m}
\]
and
\[
\Lambda_{0,n}(e_{m_1},\ldots,e_{m_n})\,=\,r p_*\cv
\]
where $p:\,\M_{g,n}^{1/r,\bm}\,\to\,\M_{g,n}$, 
$\bm:=(m_1,\ldots,m_n)$, and $0\leq m_i \leq r-2$. The 
corresponding potential function $\Phi_0(\bx)$ (of the Frobenius 
manifold) is given by the formula 
\[
\Phi_0(\bx)\,:=\,\sum_{n\,\geq\,3}\,\frac{1}{n!}\,
x^{m_1} \ldots x^{m_n}\,r\,\int_{\M_{0,n}^{1/r,\bm}}\,\cv,
\]
where $\bx\,:=\,(x^0,\ldots,x^{r-2})$ are coordinates 
corresponding to the basis $\{\,e_0,\ldots,e_{r-2}\,\},$ and 
$\bm\,:=\,(m_1,\ldots,m_n)$.  Summation over $m_i=0,\ldots,r-2$ 
is assumed, and the integral is understood in the sense of 
orbifolds. 

Furthermore, the potential $\Phi_0(\bx)$ is equal to the large 
phase space potential $\Phi_0(\bt)$ from Theorem \ref{thm:gwconj} 
when $t_a^m$ is set to zero for $a > 0$ and $x^m := t_0^m$ for 
all $m\,=\,0\,\ldots,r-2$. The potential $\Phi_0(\bt)$ satisfies 
the $g=0$ topological recursion relations 
\[
\frac{\partial^3 \Phi_0}{\partial t_{a_1+1}^{m_1} \partial t_{a_2}^{m_2}
  t_{a_3}^{m_3}}\,=\, \sum_{m_+,m_-}\,\frac{\partial^2 \Phi_0}{\partial
  t_{a_1}^{m_1} \partial t_{0}^{m_+}}\, \eta^{m_+,m_-}\,\frac{\partial^3
  \Phi_0}{\partial t_0^{m_-}\,\partial t_{a_2}^{m_2}\,\partial t_{a_3}^{m_3}}.
\]
The potential also satisfies the string equation, equation
(\ref{eq:string}), where the symbol $\Phit(\bt)$ is replaced by
$\Phi_0(\bt)$.
\end{thm}

This result was conjectured by Witten \cite{W,W2} before the 
proper moduli spaces had even been constructed. The proof of the 
$g=0$ \cft\ property uses intersection theory on $\M_{0,n}^{1/r}$ 
which is rather involved because of its orbifold (or stacky) 
nature. The proof of the topological recursion relations follows 
from the properties of the classes involved when they are 
restricted to boundary strata and from a presentation of the 
$\psi$ classes in terms of boundary classes. The string equation 
follows from the lifting properties of the classes under the map 
of ``forgetting'' a point labeled by $m=0$ and from the grading. 
Finally, one obtains the structure of a Frobenius manifold on 
$(\chr,\eta)$ because $\Phi_0(\bx)$ is a polynomial. We refer the 
interested reader to \cite{JKV} for details.  

\begin{rem}
Strictly speaking, the state space of the $\KdV_r$ theory should 
be $(\cht,\tilde{\eta})$, where $\cht$ is an $r$-dimensional 
vector space with basis $\{\,e_0,\,\ldots,\,e_{r-1}\,\}$, and  
with a metric in this basis given by 
$\tilde{\eta}_{m_1,m_2}\,:=\,1$ if $m_1\,+\,m_2\,\equiv 
\,(r-2)\,\mod\,r$ and $0$ otherwise. However, there is an
obvious orthogonal decomposition of vector spaces 
$\cht\,=\,\chr\, 
\oplus\,\ch'$ where $\ch'$ is the one dimensional vector space 
with basis $\{e_{r-1}\}$, and because $\cv$ vanishes on 
$\M_{g,n}^{1/r,\bm}$ whenever one of the  $m_i$ is $r-1$. 
Therefore, the decomposition $\cht\,=\,\chr\, 
\oplus\,\ch'$ is a direct sum of \cfts\ where $\ch'$ is the 
trivial, one-dimensional \cft. For this reason, we can and will 
restrict ourselves to the state space $(\chr,\eta).$ 
\end{rem}

What is interesting about this approach is that it offers the possibility
of generalization to higher genera. In \cite{JKV},  axioms
were formulated for the virtual class $\cv$ in order to obtain a  \cft\ in
all genera drawing upon an analogy with Gromov-Witten invariants. These
axioms should be viewed as an analog  of those of Behrend-Manin \cite{BMa}
for the virtual fundamental class in the case of the moduli space of stable
maps.  

Finally, it is also worth observing that the Frobenius structure here
\textsl{cannot} arise as the  quantum cohomology of a smooth, projective
variety because it would correspond to fractional dimensional 
cohomology classes on the target variety.

\section{The tensor product of Frobenius manifolds}

The category of cohomological field theories has a natural tensor product
\cite{KMK} described as follows.

\begin{df}
Let $(\ch',\eta',\Lambda')$ and $(\ch'',\eta'',\Lambda'')$ be \cfts. Their
\textsl{tensor product} is
$(\ch'\,\otimes\,\ch'',\eta'\,\otimes\,\eta'',\Lambda)$ where 
\[
\Lambda_{g,n}(v'_1\otimes v''_1,\ldots, v'_n\otimes v''_n)\,:=\,
\Lambda_{g,n}'(v'_1,\ldots v'_n)\,\cup\,\Lambda''_{g,n}
(v''_1,\ldots,v''_n)
\]
for all $v'$ in $\ch'$ and $v''$ in $\ch''$.
\end{df}

This is nothing more than the fact that the diagonal map
$\M_{g,n}\,\to\,\M_{g,n}\,\times\, \M_{g,n}$ is a coproduct with 
respect to the composition maps of the modular operad 
$\{\,H_\bullet(\M_{g,n})\,\}$. We will only discuss the situation 
where $g\,=\,0$ and the potential functions are polynomial in the 
flat coordinates. 

The tensor product operation, when written in terms of the 
underlying potential functions, is highly nontrivial 
\cite{Ka,Ka2}, even in genus zero. In the case of Gromov-Witten 
invariants, Behrend \cite{B} proved that the \cft\ arising from 
$\M_{g,n}(V'\,\times\,V'')$ is the tensor product of that arising 
from $\M_{g,n}(V')$ and $\M_{g,n}(V'')$. When restricted to genus 
zero, one can view this result as a deformation of the K\"unneth 
theorem. 

There exists a notion of the tensor product of (nonformal) 
Frobenius manifolds \cite{Ka}. Because the potential functions of 
$\KdV_r$ are polynomial, their tensor product is also a (nonformal) Frobenius
manifold.  There is, however, an unanswered question. Is there a natural
moduli space whose intersection numbers assemble into a 
generating function related to the tensor products of the 
Frobenius manifolds associated to $\KdV_{r^i}$ for 
$i\,=\,1\,\ldots\,k$? These moduli spaces would be analogs of 
$\M_{0,n}(V_1\times\cdots\times V_k)$ in the theory of 
Gromov-Witten invariants.  The answer turns out to be \emph{yes}, 
as explained in the following theorem.  This theorem is an 
immediate consequence of Lemma~\ref{pb-prod} 

\begin{thm}
Let $(\ch^{(r^i)},\eta^{(i)})$ denote the Frobenius manifold 
associated to the $r^i$-th Gelfand-Dickey hierarchy, 
$\KdV_{r^i}$, for $r^i\,\geq\,2$ and $i\,=\,1\ldots,k$. Let 
$(\ch,\eta)$ denote the tensor product of these Frobenius 
manifolds.  The tensor product Frobenius structure arises from 
$\Lambda\,:=\,\{\,\Lambda_{0,n}\,\}$ where the forms 
$\Lambda_{0,n}  \in 
H^\bullet(\M_{g,n})\,\otimes\,T^n\,\ch^*$, are defined 
by 
\[
\Lambda_{0,n}(e_{m_1^1}\otimes\ldots\otimes e_{m_1^k},\ldots,
e_{m_n^1}\otimes\ldots\otimes e_{m_n^k})\,=\, (\prod_{i=1}^k r^i)
p_*c^{\bo/\br}. 
\]
and where $p\,:\,\M_{g,n}^{\bo/\br,\bbm}\,\to\,\M_{g,n}$,
$\br\,:=\,(r^1,\ldots,r^k)$, $\bbm\,:=\,(\bm^1,\ldots,\bm^k)$,
$m^i_j\,=\,0 , \ldots, r^i-2$ for all $i\,=\,1,\ldots,k$ and
$j\,=\,1,\ldots,n$.
\end{thm}

One may regard the above result as a geometric realization (a 
kind of ``A-model'') for the tensor product of these Frobenius 
manifolds. Isomorphic Frobenius manifolds can be constructed from 
versal deformations of $A_{r^i}$ singularities in \cite{Ma,Ma2} 
(a Landau-Ginzburg model). The latter construction can be 
generalized to yield a Frobenius manifold associated to any 
singularity of type A-D-E, or even more generally, to Coxeter 
groups. It is known that as Frobenius manifolds, the tensor 
product of the Frobenius manifolds associated to $A_{2}$ 
($\KdV_3$) and $A_{3}$ ($\KdV_4$) is isomorphic to the Frobenius 
manifold associated to $E_6$, while the tensor product of the  
Frobenius manifolds associated to $A_2$ ($\KdV_3$) and  $A_4$ 
($\KdV_5$) is isomorphic to the Frobenius manifold associated to 
$E_8$ \cite{DiVV,Du2}.

\bibliographystyle{amsplain}

\providecommand{\bysame}{\leavevmode\hbox to3em{\hrulefill}\thinspace}

\end{document}